\documentclass[12pt]{article}
\usepackage{a4wide}
\usepackage{amsmath}
\usepackage{amssymb}
\usepackage{booktabs}
\usepackage{latexsym}
\usepackage{mathpazo}
\usepackage{mathrsfs}
\usepackage{tikz}
\usepackage[normalem]{ulem}
\usepackage{soul}

\usepackage{todonotes}

\newcommand{\E}{{\rm \bf E}}

\newcommand{\cav}{{\rm cav\,}}

\newcommand{\rmd}{{\rm d}}

\newcommand{\calF}{{\cal F}}

\newcommand{\calS}{{\cal S}}
\newcommand{\calT}{{\cal T}}

\newcommand{\ep}{\varepsilon}

\usepackage{graphicx}
\usepackage{graphics}

\newtheorem{theorem}{Theorem}[section]

\newtheorem{definition}[theorem]{Definition}

\newtheorem{comment}[theorem]{Remark}

\newcommand{\dN}{\mathbb{N}}

\newcommand{\dR}{\mathbb{R}}
\newcommand{\prob}{\mathbb{P}}

\newcounter{figurecounter}
\begin{document}
\title{The Maxmin Value of Repeated Games with Incomplete Information on One Side and Tail-Measurable Payoffs%
\thanks{
Castellon Koltun acknowledges the support of the Israel Science Foundation grant \#722/18.
Lehrer acknowledges the support of the Israel Science Foundation grant \#591/21 and the DFG grant KA 5609/1-1.
Solan acknowledges the support of the Israel Science Foundation grant \#211/22.
The authors thank the reviewers for comments that improved the presentation.}}

\author{Gil Bar Castellon Koltun,\footnote{School of Mathematical Sciences, Tel Aviv University, Tel Aviv, Israel, 6997800, E-mail: gilbarkoltun@tauex.tau.ac.il.}\,\,
Ehud Lehrer,\thanks{Department of Economics, Durham University, Durham, DH1 3LB, UK,  E-mail: ehud.m.lehrer@durham.ac.uk.}\,\,
and
Eilon Solan\footnote{School of Mathematical Sciences, Tel Aviv University, Tel Aviv, Israel, 6997800, E-mail: eilons@tauex.tau.ac.il.}}

\maketitle

\begin{abstract}
We study two-player zero-sum repeated games with incomplete information on one side, where the payoff function is tail measurable (and not necessarily the long-run average payoff). 
We show that the maxmin value equals the concavification of the value function of the non-revealing game. 
In addition, we provide an example demonstrating that, under tail-measurable payoffs, the value of the game may fail to exist.
\end{abstract}

%\begin{abstract} \blue{I would start with the model we consider and then describe the main results.}
%We prove that for two-player zero-sum repeated games with incomplete information on one side,
%where the payoff function is tail measurable (and not necessarily the long-run average payoff),
%the maxmin value
%is equal to the concavification of the value function of the non-revealing game.
%We provide an example which demonstrates that the value need not exist when the payoff function is tail measurable.
%\end{abstract}

\noindent
\textbf{Keywords:} Repeated games with incomplete information on one side, tail-measurable payoff function, maxmin value, concavification.
\bigskip

\noindent
\textbf{JEL Classification:}
C72, C73.

\noindent
\textbf{MSC 2020 Classification:}
91A05, 91A20, 91A27.

%%%%%%%%%%%%%%%%%%%%%%%%%%%
\section{Introduction}

Since the influential work of Aumann and Maschler (1968, 1995; below AM) on repeated games with incomplete information, the optimal use of information in repeated interactions has been extensively studied.
Aumann and Maschler focused on two-player zero-sum repeated games with incomplete information on one side and the long-run average payoff.
They proved
the existence of the value, and characterized it by means of the value function of the non-revealing game,
which is the game in which neither player is informed of the state of nature.
Specifically,
they showed that the value of the zero-sum repeated game with incomplete information on one side at the initial prior $p$
is the concavification of the value function of the non-revealing game at $p$.
Moreover, 
they proved that this value is also the limit of the values of the finite-horizon games as the horizon increases to infinity.
When this property holds, the limit is also called the \emph{asymptotic value}.
Further, the $\ep$-optimal strategies that AM constructed are $2\ep$-optimal in all sufficiently long games
and in all discounted games, provided the players are sufficiently patient.
Such strategies are called \emph{uniformly $\ep$-optimal}.

Aumann and Maschler (1995, Chapter 2, postscript c) advocated the study of games with the long-run average payoff on several grounds.
These games can model situations where the game has in fact a finite horizon, but the players do not know the length of the game,
or they do know the length of the game but taking it into account in forming the strategy is too cumbersome or time consuming. 
Additionally, 
the existence of uniform $\ep$-optimal strategies allows us to focus on the optimal behavior of the players.
Put differently, 
games with the long-run average payoff allow us to study rules of thumb -- strategies that are approximately optimal for all finite-horizon versions of the game that are sufficiently long.

The economic literature typically studies models in which players receive stage payoffs
that are aggregated into a utility function, for example, by taking their discounted sum,
their average over a finite horizon, or their long-run average.
In some settings, however, players’ objectives cannot be captured by such evaluations.
For instance, an investor may care not only about the long-run average profit per period
but also about the average stage-to-stage fluctuations,
and a player may care about her long-run average payoff
as well as her performance relative to other players.
Similarly, in models of R\&D, the outcome may depend on whether the frequency of positive breakthroughs exceeds a certain threshold, or on whether negative results do not cluster too often.

The work of AM
%Aumann and Maschler (1968, 1995)
was influential on the literature that studies information revelation,
such as Bayesian persuasion and cheap talk (see, e.g., Kamenica and Gentzkow (2011) and Forges (2020)).
It is therefore important to extend AM's results to zero-sum repeated games with incomplete information on one side where the payoff function is more general than the long-run average payoff.

As exhibited by, e.g., Sorin (1984, 1985) and Laraki (2000), 
repeated games with incomplete information on one side and a general payoff function may fail to have a value;
that is, the minmax value may differ from the maxmin value.%
\footnote{Though the models studied by Sorin (1984, 1985) and Laraki (2000) allow the state to change dynamically along the game,
since transitions in these models are deterministic,
the models
can be recast as repeated games with a state that remains constant along the game and a general payoff function.}
A natural question is, then, whether the existence of the value and its characterization extend to other classes of 
repeated games with incomplete information on one side with economically relevant payoff functions.
Put differently, 
we would like to know which properties of the payoff function ensure the existence of the value
and its characterization as the concavification of the value function of the non-revealing game.

\paragraph{The paper's contribution.}

We study two-player zero-sum repeated games with incomplete information, 
where the state of nature is realized at the outset and remains fixed throughout the game.
Only Player~I observes the realized state. 
The payoff that Player~II pays Player~I is determined at the end of the game and depends on the entire infinite history in a manner that makes it tail measurable.
Those are payoff functions that depend on the long-run behavior of the players, and are invariant under changes to any finite prefix of the play.
Both the long-run average payoff
and the examples mentioned two paragraphs above are tail measurable, but the discounted payoff is not.

%In this paper we study repeated games with incomplete information on one side, 
%where the payoff function is tail measurable.
%Those are payoff functions that depend on the long run behavior of the players, and are invariant under changes to any finite prefix of the play.
%The long-run average payoff 
%as well as the two examples mentioned two paragraphs above are tail measurable, but the discounted payoff is not.

We show that 
%for this class of payoff functions,
%the value may fail to exist.
%However, 
the maxmin value,
namely, the amount that Player~I
%, the maximizer who is informed of the state of nature,
can guarantee,
%Aumann and Maschler's (1968, 1995) characterization:
is equal to the concavification of the value function of the non-revealing game.
%By virtue of an example, we also show that when the payoff function is tail measurable,
%the value might not exist.
%
%One consequence of this result concerns the information revelation process required from Player~I.
%
This implies that to guarantee the maxmin value,
Player~I needs to transmit information to Player~II only once, at the beginning of the game.
After that information revelation,
Player~I does not use or transmit additional information.
In particular, the maxmin strategy of Player~I has a particularly simple form.

%\blue{I would skip this paragraph, as it is not a consequence and it is elaborated upon later.} Second, 
%in our proof we construct an $\ep$-optimal response of Player~II to any given strategy of Player~I.
%Under this construction, 
%Player~II plays in blocks;
%in each block Player~II plays as if no information is going to be revealed in the future,
%and, if a significant amount of information is revealed,
%then Player~II updates her belief about the state of nature and starts a new block, where she follows a new strategy that is good against the new belief, assuming no information will be revealed in the future.
%As we elaborate below, 
%this suggests an information-based good response for Player~II,
%rather than a payoff-based response as in AM.

We also show that imposing a natural condition such as tail measurability,%
\footnote{In fact, the payoff function in our example where the value does not exist satisfies a stronger condition than tail measurability, called \emph{shift invariance}. The payoff function is \emph{shift invariant} if it is invariant to shortening the play by removing a finite number of action pairs.}
is not sufficient to guarantee the existence of the value in two-player zero-sum games with incomplete information on one side. 
Specifically, we construct an alternating-move two-player zero-sum repeated game with incomplete information on one side, which is a variant of a game studied by Laraki (2000), in which the payoff is tail measurable but the game does not possess a value.

\paragraph{The idea of the proof.}

%Aumann and Maschler's (1968, 1995) 
AM's
arguments show that the maxmin value is always at least the concavification of the value function of the non-revealing game.
To prove the main result, we need to prove the reverse inequality.

To this end, we fix a strategy of Player~I,
and construct a response of Player~II that guarantees that the expected payoff is not much higher than the concavification of the value function of the non-revealing game.

We consider the sequence of posterior beliefs of Player~II, which are determined by the strategy of Player~I and the play,
and we define Player~II's response in blocks:
in each block, Player~II follows an $\ep$-optimal strategy
in the non-revealing game,%
\footnote{
In fact, Player~II will follow 
a \emph{subgame $\ep$-optimal} strategy; 
namely, 
a strategy that guarantees the value of the non-revealing game up to $\ep$ in all subgames.}
whose prior is equal to Player~II's posterior belief at the beginning of the block,
and the block ends as soon as the posterior belief 
deviates 
significantly from its value at the beginning of the block.
%considerably from the posterior belief at the beginning of the block.
Since the posterior belief is a martingale, it converges a.s.,
and hence the number of blocks is a.s.~finite.%
\footnote{
In the model with the long-run average payoff,
a similar construction, with blocks of size $1$,
was used in Sorin (2002, Section~3.4.1)
to bound the difference between the values of the finite-horizon games and the concavification of the value function of the non-revealing game.}
The reason for ending a block when the posterior belief deviates considerably from the posterior belief at the beginning of the block and re-initiating the strategy of Player~II,
is that Player~II's strategy within a block is $\ep$-optimal in the non-revealing game with prior that is equal to the posterior belief at the beginning of the block.
When the posterior belief deviates considerably from this initial prior,
the strategy that Player~II follows is no longer good, and hence should be revised.

We then show that when the block is infinite with high probability,
the expected payoff is not much larger than the value of the non-revealing game 
with prior which is the belief at the beginning of the block.
This part is not trivial,
because Player~II follows an $\ep$-optimal strategy in the non-revealing game,
while Player~I plays a strategy that does depend on the state of nature,
albeit in a way that does not significantly change Player~II's belief.

Using Jensen's inequality,
we deduce that the expected payoff in the whole game is not much larger than the concavification of the value function of the non-revealing game,
evaluated at the initial prior.

\paragraph{Comparison of our maxmin strategy and that of Aumann and Maschler (1968, 1995).}

There is an important conceptual difference between the approximate minmax strategy that we construct for Player~II and the one constructed by 
%Aumann and Maschler (1968, 1995) 
AM
for the long-run average payoff.
The strategy constructed by 
AM
%Aumann and Maschler (1968, 1995)
is \emph{payoff-based}:
it requires Player~II to adapt her mixed action at any given stage to the vector of payoffs she could obtain in the last stage
(where each coordinate of the vector corresponds to one state).
When the payoff function is general,
there are no stage payoffs,
and hence the strategy we construct is \emph{information-based}:
Player~II adapts her play only when the amount of information revealed by Player~I is sufficiently large.
To be able to calculate her posterior belief,
Player~II needs to know Player~I's strategy,
and hence the strategy we construct defends the maxmin value, and is not proper without the knowledge of Player~I's strategy -- there does not necessarily exist a single strategy of Player~II that is good against \emph{every} strategy of Player~I.

\paragraph{The role of tail measurability of the payoff.}

We would like to highlight that the tail measurability of the payoff function
affects the argument for the two players in different ways.

With regard to the informed Player~I,
tail measurability implies that the information she wishes to reveal along the game 
can be equivalently revealed at the outset of the game,
and there is no reason to reveal it piecemeal. 
Indeed, 
if Player~I revealed her information with some delay,
Player~II could wait until all information is revealed, play optimally afterwards,
and the tail measurability of the payoff function implies that this delay does not affect the total payoff.

With regard to the uninformed Player~II,
tail measurability implies that the strategy described above which plays in blocks is good:
the number of blocks is finite a.s.,
in infinite blocks the strategy yields low payoff,
and since payoffs are tail measurable, 
the blocks with finite length do not affect the payoff that Player~II can guarantee in the game. 

\paragraph{Related literature.}

Our paper is related to several strands of literature.

When at each stage the outcome is a stage payoff,
it is natural to study the asymptotic behavior of the value of the finite horizon games.
This includes the convergence of the value of the finite horizon games as the length of the game increases,
the relation between the limit -- the asymptotic value -- and the value of the infinite horizon game, the examination of the rate of convergence to the limit,
and the existence of uniform $\ep$-optimal strategies that are good in all sufficiently long games. 
Two other research directions are the study of games
where both players are uninformed
and games where the state changes along the play.  

%Our paper is related to the rich literature that studies when the asymptotic value exists,
%and whether players have uniformly $\ep$-optimal strategies.
A few influential papers in these directions include
Mertens and Zamir (1971, 1980), 
who proved that for repeated games with incomplete information on both sides the asymptotic value need not exist, and characterized the maxmin and the minmax values;
Sorin (1984, 1985),
who did the same for two-player zero-sum Big Match games with incomplete information on one side;
Rosenberg and Vieille (2000),
who proved that the asymptotic value exists in recursive games;
Gensbittel, Oliu-Barton, and Venel (2014),
who proved the existence of the asymptotic value in stochastic games with a more informed controller;
Renault (2006),
who showed that in Markov games the asymptotic value exists and both players have uniform $\ep$-optimal strategies;
Lehrer and Shaiderman (2022), who studied games where the state evolves according to a Kronecker system,
and Ziliotto (2024), who showed that the asymptotic value exists in absorbing games.
On the negative side,
Sorin and Zamir (1985)
showed that the asymptotic value may fail to exist when Player~I knows the state of nature
but does not know the belief of Player~II on the state of nature;
and Ziliotto (2016) showed that the asymptotic value may fail to exist when the state variable is partially observed.
For surveys on this topic, see
Sorin (2003) and Mertens, Sorin, and Zamir (2015).

The paper also adds to the literature that studies the value of infinite-horizon zero-sum games with general payoff functions.
This class of games was first studied by Gale and Stuart (1953),
who proved that in alternating-move games with a closed or open winning set,
one of the players has a winning strategy.
Blackwell (1969) studied the analogous model when players move simultaneously,
and proved the existence of the value when the winning set is $G_\delta$,
namely, a countable intersection of open sets.
Martin (1975, 1998) extended these two results to the case
where the outcome is any bounded and Borel measurable function.
Arieli and Levy (2011) proved that two-player zero-sum alternating-move games with incomplete
information on one side,
where the payoff function is the characteristic function of some set, need not have a value.
They also showed that in their model with two states of nature, the game does admit a value.

The advantage of tail-measurable 
%and shift-invariant 
payoff functions in the study of repeated games with complete information has been recently exhibited in several papers.
Ashkenazi-Golan, Flesch, Predtetchinski, and Solan (2022) proved the existence of $\ep$-equilibria
in repeated games with countably many players and tail-measurable payoff functions.
Ashkenazi-Golan, Flesch, Predtetchinski, and Solan (2025) provided a Folk Theorem for repeated games with finitely many players and tail-measurable payoff functions.
And Flesch and Solan (2023) proved that two-player non-zero-sum stochastic games with shift-invariant payoffs admit $\ep$-equilibria.

\paragraph{Structure of the paper.}

The paper is organized as follows.
The model and the main result are presented in Section~\ref{section:model}.
The proof that the maxmin value coincides with the concavification of the value function of the non-revealing game
is given in Section~\ref{section:proof}.
Two examples that demonstrate that two-player zero-sum games with incomplete information on one side and tail-measurable payoffs
need not have a value are provided in Section~\ref{section:examples}.
Section~\ref{section:discussion} lists some open problems.

%%%%%%%%%%%%%%%%%%%%%%%%%%%
\section{The Model}
\label{section:model}

\begin{definition}
A \emph{two-player zero-sum repeated game with incomplete information on one side} is a $5$-tuple
$(K,I,J,p,f)$,
where $K$ is a finite set of states of nature,
$I$ and $J$ are finite sets of actions for the two players,
$p \in \Delta(K)$ is the prior distribution,%
\footnote{For every finite set $X$, $\Delta(X)$ is the set of probability distributions on $X$.}
and $f : K \times (I \times J)^\dN \to \dR$ is the payoff function,%
\footnote{Throughout the paper, $\dN = \{0,1,2,\ldots\}$ and $(I \times J)^0$ is the singleton containing the empty history.}
which is assumed to be bounded and Borel measurable.
%We will denote the game by $\Gamma(p)$.
\end{definition}

The game is played as follows.
Nature selects a state of nature $k^* \in K$ according to $p$,
and reveals $k^*$ to Player~I, but not to Player~II.
At every stage $t \in \dN$,
the two players, after observing their past actions, simultaneously select actions $i_t\in I$ and $j_t \in J$.
The payoff that Player~I tries to maximize and Player~II tries to minimize is $f(k^*,i_0,j_0,i_1,j_1,\ldots)$.

The players have perfect recall and full observation, hence a \emph{strategy} of Player~I is a function 
$\sigma : K \times  \left(\bigcup_{t \in \dN} (I \times J)^{t}\right) \to \Delta(I)$
%$\sigma : \bigcup_{t \in \dN} \bigl(K \times (I \times J)^{t}\bigr) \to \Delta(I)$
and a \emph{strategy} of Player~II is a function
$\tau : \bigcup_{t \in \dN} (I \times J)^{t} \to \Delta(J)$.
We denote by $\calS$ and $\calT$ the sets of strategies of the two players, respectively.
Every pair $(\sigma,\tau) \in \calS \times \calT$,
together with the prior $p$,
induces a probability distribution $\prob_{p,\sigma,\tau}$ on the set of \emph{plays} $H^\infty := K \times (I \times J)^\dN$
(equipped with the product sigma-algebra, denoted $\calF$).
Expectation w.r.t.~this probability distribution is denoted $\E_{p,\sigma,\tau}[\,\cdot\,]$.

%As the analysis of Aumann and Maschler (1968, 1995) reveals,
%the study of the game $\Gamma$ requires the study of the same game with all possible priors.
From now on we fix $K$, $I$, $J$, and $f$,
and for every $p \in \Delta(K)$ we denote by $\Gamma(p)$ the two-player zero-sum repeated game with incomplete information on one side with prior $p$.
We will assume w.l.o.g.~that $0 \leq f \leq 1$.

The \emph{maxmin value} of $\Gamma(p)$ is
\begin{equation}
\label{equ:maxmin}
\underline v(p) = \sup_{\sigma\in \calS}\inf_{\tau \in \calT} \E_{p,\sigma,\tau}[f],
\end{equation}
and the \emph{minmax value} of $\Gamma(p)$ is
\begin{equation}
\label{equ:minmax}
\overline v(p) = \inf_{\tau \in \calT} \sup_{\sigma \in \calS}\E_{p,\sigma,\tau}[f].
\end{equation}
The game $\Gamma(p)$ admits a \emph{value}, denoted $v(\Gamma(p))$, if its minmax and maxmin values coincide.
A strategy $\sigma$ (resp., $\tau$) that attains the supremum (resp., infimum) in Eq.~\eqref{equ:maxmin} (resp., Eq.~\eqref{equ:minmax}) up to $\ep$
is called an \emph{$\ep$-maxmin} (resp., \emph{$\ep$-minmax}) strategy.

\begin{comment}
When the payoff function is the \emph{long-run average of stage payoffs}, namely,

\[ f(k,i_0,j_0,i_1,j_1,\ldots) = \limsup_{T \to \infty}\frac{1}{T} \sum_{t=0}^{T-1} g_k(i_t,j_t), \]
for some collection of functions $g_k : I \times J \to \dR$, $k \in K$,
the literature defines the maxmin value (resp.,~minmax value) as the amount that Player~I (resp., Player~II) can guarantee in all sufficiently long games
and Player~II (resp., Player~I) can defend in all sufficiently long games, see, e.g., Rosenberg and Vieille (2000).
With this definition,
the existence of these two values is not guaranteed and requires a proof.
However, 
by AM
they coincide with those given in Eqs.~\eqref{equ:maxmin} and~\eqref{equ:minmax}.

When the payoff function is general, as in our case,
it is not clear how to define the finite-stage game,
and the definitions we provided to the maxmin and minmax values are natural
and always exist.
\end{comment}

\begin{definition}
Let $\Gamma(p)$ be a two-player zero-sum repeated game with incomplete information on one side.
The \emph{non-revealing game} with prior $p$, denoted $\Gamma^{\rm NR}(p)$,
is the game that is similar to $\Gamma(p)$,
except that both players have incomplete information:
no player is told the identity of the state of nature $k^*$ that is selected at the outset of the game.
As a result, a strategy of Player~I is a function
$\sigma : \bigcup_{t \in \dN} (I \times J)^{t} \to \Delta(I)$.
\end{definition}

The game $\Gamma^{\rm NR}(p)$ is a two-player zero-sum repeated game (with complete information),
where
the action sets of the two players are $I$ and $J$, respectively, and
the payoff function is $\sum_{k \in K} p(k)f(k,\cdot)$.
By Martin (1998),
this game admits a value, denoted $u(p)$.
Since $0 \leq f \leq 1$, the function $p \mapsto u(p)$ is continuous and even $1$-Lipschitz,
see Part~2 of the proof of Theorem~\ref{theorem:main}.

%\begin{lemma}
%\label{lemma:symmetric:value:exist}
%For every prior $p \in \Delta(K)$ and bounded and Borel measurable payoff function $f$,
%the game $\Gamma^{\rm NR}(p)$ admits a value, denoted $u(p)$.
%\end{lemma}

Denote by $\cav u$ the \emph{concavification} of $u$,
namely, the least real-valued concave function on $\Delta(K)$ that is larger than or equal to $u$.
AM
%Aumann and Maschler (1968, 1995)
proved that when $f$ is the long-run average payoff,
the value of the game $\Gamma(p)$ exists and is equal to $(\cav u)(p)$.
Their argument does in fact show that $\underline v(p) \geq (\cav u)(p)$
for every bounded and Borel measurable payoff function $f$.

\begin{definition}
The payoff function $f$ is \emph{tail measurable} if
\[ f(k,i_0,j_0,i_1,j_1,\dots) = f(k,i'_0,j'_0,i'_1,j'_1,\dots), \]
for every $k \in K$ and every two sequences $(i_0,j_0,i_1,j_1,\dots), (i'_0,j'_0,i'_1,j'_1,\dots) \in (I \times J)^\dN$
that differ only in finitely many coordinates.
\end{definition}

The limsup and liminf of the average stage payoffs,
as well as the limsup and liminf of stage payoffs, are examples of tail-measurable payoff functions.
Criteria that involve convergence to a target set,
like in Blackwell approachability (Blackwell, 1956), 
are also tail measurable.
%Two additional instances where the payoff function is tail measurable are when players have several objectives, each is tail measurable,
%and they either would like to satisfy all of them
%(Basset et al., 2015),
%or have a lexicographic ranking among the objectives (Chatterjee et al., 2024).
%The discounted payoff is not tail measurable.
Tail measurability also arises when players face several tail-measurable objectives and either require all of them to be satisfied (Basset et al., 2015) or compare them lexicographically (Chatterjee et al., 2024).
%Two additional instances in which the payoff function is tail measurable occur when players have multiple objectives, each of which is tail measurable. 
%In addition, the players either aim to satisfy all objectives (Basset et al., 2015),
%or follow a lexicographic ranking among them (Chatterjee et al., 2024). 
Finally, various classical winning conditions in the computer science literature, such as the B\"uchi, co-B\"uchi, and parity
%, Streett, and M\"uller
(see, e.g., Chatterjee and Henzinger, 2012), are also tail measurable.%
\footnote{
We say that Player~I wins according to the B\"uchi (resp., co-B\"uchi) winning condition if some (resp., no) action profile in a given set of action profiles is selected infinitely often.
Player~I wins under the parity condition if, for a fixed ranking (priority function) on action profiles, the least rank that occurs infinitely often along the play is even. 
}
In contrast, the discounted payoff is not tail measurable.

Our main result states that when $f$ is tail measurable, the maxmin value of $\Gamma(p)$ is equal to $(\cav u)(p)$.

\begin{theorem}
\label{theorem:main}
If $f$ is tail measurable,
then $\underline v(p) = (\cav u)(p)$, for every $p \in \Delta(K)$.
\end{theorem}

The tail measurability in 
Theorem~\ref{theorem:main} cannot be discarded.
Indeed,
Sorin (1984, 1985) provided an example where $f$ is not tail measurable,
the value does not exist, and the maxmin value is not equal to $\cav u$.
In Section~\ref{section:examples} we provide examples of two-player zero-sum repeated games with incomplete information on one side
where $f$ is shift invariant (and so, in particular, tail measurable),
and for which the value does not exist.

%%%%%%%%%%%%%%%%%%%%%%%%%%%%%%%%%%%%%%%%%%%%%%%%%%%%%%%%%%%%%%%%%%%%%%%%
\section{Proof of Theorem~\ref{theorem:main}}
\label{section:proof}

The proof is divided into 12 parts.
The first handles the inequality $\underline v(p) \geq (\cav u)(p)$,
whose proof is similar to the analogous inequality in  
AM.
%Aumann and Maschler (1968, 1995)
The remaining parts handle
the more challenging inequality $\underline v(p) \leq (\cav u)(p)$.

\bigskip
\noindent\textbf{Part~1: $\underline v(p) \geq (\cav u)(p)$.}

By never using her information, the informed player can guarantee $u(p)$ in the game $\Gamma(p)$.
Hence $\underline v(p) \geq u(p)$.
%\st{Aumann and Maschler} 
AM's splitting lemma implies that the function $\underline v$ is concave in $p$
(an alternative argument 
for the concavity of $\underline v$
is provided in Proposition 2.2 in Sorin (2002)).
Since $\cav u$ is the smallest concave function that is greater than or equal to $u$,
it follows that $\underline v(p) \geq (\cav u)(p)$.

\bigskip

We now turn to the proof that $\underline{v}(p) \leq (\cav u)(p)$. 
As outlined in the introduction, the key idea is as follows. 
We fix a strategy $\sigma$ for the informed player and, for every $\ep > 0$,
we will construct a corresponding response $\tau^*$ for the uninformed player
such that $\E_{p,\sigma,\tau^*}[f] \leq (\cav u)(p) + 8|K|\ep$. 
Since $\ep$ and $\sigma$ are arbitrary, this will imply that $\underline{v}(p) \leq (\cav u)(p)$. 

To construct $\tau^*$, we will partition time into random-length blocks. 
Within each block, the strategy will follow a subgame $\ep$-optimal
%minmax 
strategy in the non-revealing game, with respect to the belief at the beginning of the block. The block will terminate when Player II's belief deviates significantly from her belief at the beginning of the block.

The proof will be by induction on the number of states of nature $|K|$.
When $K=\{k^*\}$,
there is no information to reveal.
Moreover, the game is a standard Blackwell game,
and hence by Martin (1998) its value exists, and it is equal to $u(\mathbf{1}_{k^*}) = (\cav u)(\mathbf{1}_{k^*})$.

We next assume by induction that the claim holds in all repeated games with incomplete information on one side with tail-measurable payoff function and fewer than $|K|$ states.

We start with a few observations, notations, and definitions.

\bigskip
\noindent\textbf{Part~2:
The function $\cav u$ is uniformly continuous.}
\newline
For every prior $p \in \Delta(K)$ and every pair of strategies $(\sigma',\tau')$ in the non-revealing game $\Gamma^{\rm{NR}}(p)$,
\[ \E_{p,\sigma',\tau'}[f]
= \sum_{k \in K} p(k) \E_{k,\sigma',\tau'}[f]. \]
Since $0 \leq f \leq 1$,
the function $p \mapsto \E_{p,\sigma',\tau'}[f]$ is 1-Lipschitz.
Hence, 
the function
\[ p \mapsto u(p) = \sup_{\sigma'} \inf_{\tau'} \E_{p,\sigma',\tau'}[f] \]
is 1-Lipschitz as well, 
and in particular continuous.
Hence,
the function $\cav u$ is continuous.%
\footnote{To see this implication, 
consider a sequence $(p_n)_{n \in \dN}$ that converges to a limit $p$.
Suppose that for every $n \in \dN$ we have $(\cav u)(p_n) = \sum_{k=1}^{|K|} \alpha_{n,k} u(p_{n,k})$,
where $(\alpha_{n,k})_{k=1}^{|K|}$ are non-negative and sum to $1$.
W.l.o.g.~assume that the limits $\alpha_k := \lim_{n \to \infty} \alpha_{n,k}$ and $p_k := \lim_{n \to \infty} p_{n,k}$ exist for every $k$.
Since $u$ is continuous,
$(\cav u)(p) \geq \sum_{k=1}^{|K|} \alpha_k u(p_k) = \lim_{n \to \infty} (\cav u)(p_n)$.
Hence, $\cav u$ is upper-semicontinuous.
Since any upper-semi-continuous and concave function 
that is defined on a compact convex set is continuous,
it follows that $\cav u$ is continuous.}
Since $\Delta(K)$ is compact, $\cav u$ is uniformly continuous.

\bigskip
\noindent\textbf{Part~3: The algebra $\calF^{\rm II}_t$ that is generated by public histories of length $t$.}

A \emph{public history} is an element of the set $H_{\rm pub} := \bigcup_{t \in \dN} (I\times J)^{t}$.
The \emph{length} of a public history $h = (i_0,j_0,i_1,j_1,\dots,i_{t-1},j_{t-1})$ is $t$.
When we need to emphasize the length of the public history, we will denote a public history of length $t$ by $h_t$.

%and an infinite public play is an element of the set $H^\infty_{\rm pub} := (I \times J)^\infty$.
We say that a public history $h' \in H_{\rm pub}$ (resp., a play $h^\infty \in H^\infty$) extends a public history $h$,
and write $h \prec h'$ (resp., $h \prec h^\infty$),
if $h$ is a strict prefix of $h'$ (resp., $h^\infty$).
The concatenation of two public histories
$h = (i_0,j_0,\ldots,i_{t-1},j_{t-1})$
and
$h' = (i'_0,j'_0,\ldots,i'_{t'-1},j'_{t'-1})$
is
$h \circ h' := (i_0,j_0,\ldots,i_{t-1},j_{t-1},i'_0,j'_0,\ldots,i'_{t'-1},j'_{t'-1})$.

For every public history $h$,
denote by $C(h)$ the set of all plays that extend $h$:
\[ C(h) := \bigl\{ h^\infty \in H^\infty \colon h \prec h^\infty\bigr\}. \]

For every $t \in \dN$, denote by $\calF^{\rm II}_t$ the algebra on $H^\infty$ generated by the sets $C(h_t)$, where $h_t$ ranges over all public histories with length $t$.
Denote by $\calF^{\rm II}$
%$ := \sigma\left( \bigcup_{t \in \dN} \calF^{\rm II}_t\right)$ 
the sub-sigma-algebra of $\calF$ generated by the family
$\{C(h) \colon h \in H_{\rm pub}\}$.

\bigskip
\noindent\textbf{Part~4: Restriction of strategies and payoff functions to subgames.}

Abusing terminology, 
for every public history $h \in H_{\rm pub}$
we will call the part of the game conditional that $h$ occurred the \emph{subgame that starts at $h$}.

For every strategy $\sigma' \in \calS$ and every public history $h \in H_{\rm pub}$,
let $\sigma'_h$ be the restriction of $\sigma'$ to the subgame that starts at $h$, namely, the strategy defined by
\[ \sigma'_h(k,h') := \sigma'(k, h \circ h'), \ \ \ \forall\ k \in K, \forall\ h' \in H_{\rm pub}. \]
Similarly, for every strategy $\tau' \in \calT$ and every public history $h \in H_{\rm pub}$,
let $\tau'_h$ be the restriction of $\tau'$ to the subgame that starts at $h$:
\[ \tau'_h(h') := \tau'(h \circ h'), \ \ \ \forall\ h' \in H_{\rm pub}. \]
For every public history $h$,
denote by $f_h$ the payoff function restricted to the subgame that starts at $h$:
\[ f_h(k,h') := f(k,h \circ h'), \ \ \ \forall\ k \in K, \forall\ h' \in H_{\rm pub}. \]

\bigskip
\noindent\textbf{Part~5: Posterior belief.}

For every probability distribution $p \in \Delta(K)$, 
%every strategy $\sigma' \in \calS$ of Player~I, 
every state of nature $k \in K$,
and every public history $h = (i_0,j_0,i_1,j_1,\ldots,i_{t-1},j_{t-1}) \in H_{\rm pub}$,
denote the probability of $(k,h)$ under $(p,\sigma)$ by
\[ \prob_{p,\sigma}(k,h) := p(k) \prod_{l=0}^{t-1} \sigma(k,i_0,j_0,\dots,i_{l-1},j_{l-1})(i_{l}). \]
By Bayes rule, the \emph{posterior belief} at $h$ under $(p,\sigma)$,
denoted $\pi(p,\sigma,h) \in \Delta(K)$, is given by
\begin{equation}
\label{equ:posterior}
\pi(p,\sigma,h)(k) := \frac{\prob_{p,\sigma}(k,h)}{\sum_{k' \in K} \prob_{p,\sigma}(k',h)}, \ \ \ \forall\ k \in K.
\end{equation}
The posterior belief is well-defined whenever the denominator in Eq.~\eqref{equ:posterior} does not vanish.
Note that 
$\pi(p,\sigma,h)$ does not depend on the actual initial state $k^*$ nor on Player~II's strategy.

We will denote by $p_t(h_t)$ the posterior belief at the public history $h_t$,
and emphasize that this random variable depends on both the initial belief $p$ and Player~I's strategy $\sigma$.
Since the posterior belief is independent of the state of nature,
we will view the process $(p_t)_{t \in \dN}$ as a stochastic process on $H^\infty$.

%It will be useful to consider the posterior belief as a stochastic process on $H^\infty$ 
%(which is independent of the state of nature).
%We denote by $p_t$ the random variable which indicates the posterior belief at the beginning of stage $t$,
%and emphasize that this random variable depends on both the initial belief $p$ and Player~I's strategy $\sigma$.

%For convenience,
%when $X_t$ is an $\calF^{\rm II}_t$-measurable random variable
%and $h$ is a public history of length at least $t$,
%we sometimes write $X_t(h)$ instead of $X_t(h^\infty)$,
%for any play $h^\infty \in H^\infty$ that extends $h$.
%In particular, $p_t(h_t)$ is the posterior belief at the public history $h_t$.

\bigskip
\noindent\textbf{Part~6: Subgame $\ep$-optimal
%minmax 
 strategies in the non-revealing game.}

For every $p' \in \Delta(K)$ and every $f' : H^\infty \to \dR$,
denote by $\Gamma^{\rm NR}(p',f')$
the non-revealing game with set of states $K$,
sets of actions $I$ and $J$, prior $p'$, and payoff function $f'$.
This game is a standard Blackwell game,
and hence, by Martin (1998), its value, denoted by $u(p',f')$, exists.
Moreover, by Mashiah-Yaakovi (2015),
Flesch, Herings, Maes, and Predtetchinski (2021),
or Flesch and Solan (2024),
Player~II has a subgame $\ep$-optimal
%minmax 
 strategy
in $\Gamma^{\rm NR}(p',f')$, for every $\ep > 0$;
that is, a strategy $\tau$ such that
for every public history $h_t \in H_{\rm pub}$,
the strategy $\tau_{h_t}$ is $\ep$-optimal
%minmax 
in the game
$\Gamma^{\rm NR}(p',f_{h_t})$.

The tail measurability of $f$ implies that the value of $\Gamma^{\rm NR}(p',f_h)$ depends only on the initial prior $p'$
and is independent of the public history $h$, see Ashkenazi-Golan, Flesch, Predtetchinski, and Solan (2022, Theorem 3.1):
\begin{equation}
\label{equ:value:independent}
u(p',f_h) = u(p',f), \ \ \ \forall\ p' \in \Delta(K), \forall\ h \in H_{\rm pub}.
\end{equation}

\bigskip
\noindent\textbf{Part~7: Beliefs close to the boundary of $\Delta(K)$.}

For every probability distribution $p \in \Delta(K)$
and every $\delta \in (0,1)$, 
the set of states of nature that are assigned a probability larger than $\delta/|K|$ under $p$ is
\[ K_\delta(p) := \{ k \in K \colon p(k) > \delta/|K|\}. \]
Since $\delta < 1$,
the set $K_\delta(p)$ is not empty.
Denote
by $\xi_\delta(p)$ the
probability distribution on $K$ defined by:
%projection of $p$ to $K_\ep(p)$:
\[ \xi_\delta(p)(k) := \left\{
\begin{array}{lll}
0, & & \hbox{if } k \not\in K_\delta(p),\\
\frac{p(k)}{\sum_{k' \in K_\delta(p)} p(k')}, & & \hbox{if } k \in K_\delta(p).
\end{array}
\right.
\]
Thus, all states $k$ that are assigned probability at most $\delta/|K|$ under $p$,
are assigned probability 0 under $\xi_\delta(p)$.
Note that if $K_\delta(p) \subsetneq K$,
then
$p$ is $\delta$-close
(in the $L_\infty$-norm) to the boundary of $\Delta(K)$.
The first stage in which the posterior belief is $\delta$-close to the boundary of $K$ is
\[ \theta := \min\{ t \in \dN \colon K_\delta(p_t) \subsetneq K\}. \]
In particular, $|K_\delta(p_\theta)| < |K|$ whenever $\theta < \infty$.
In the response we will construct for Player~II,
at stage $\theta$ Player~II will start following a strategy that guarantees $(\cav u)(\xi_\delta(p_\theta(h_\theta)))$ in the subgame $\Gamma(\xi_\delta(p_\theta(h_\theta)),f_{h_\theta})$,
where $h_\theta$ is the realized public history at stage $\theta$.
%which exists by the induction hypothesis.

%The probability distribution $p$ is \emph{$\ep$-good} if $p = \xi_\delta(p)$;
%that is, if there is no $k \in K$ such that $0 < p(k) < \ep$.
%Otherwise, $p$ is \emph{$\ep$-bad}.

\bigskip
\noindent\textbf{Part~8: Definition of a sequence of the random blocks.}

We now define an increasing sequence of stopping times $(t_n)_{n}$ that divides the play 
before stage $\theta$
into
(finitely or countably infinitely many)
blocks of random size.
The stopping times are defined in such a way that the process $(p_t)_{t \in \dN}$ is almost constant within each block.
%The strategy $\tau^*$ that performs well against $\sigma$ will be defined using these data.

%As we will see, $p^a_{t_n}$,
%the value of the approximating posterior at the beginning of a block,
% will be either $\ep$-good or undefined.
%We will show that the event that it is undefined occurs with low probability,
%and hence does not significantly affect the payoff.
Fix $\ep \in (0,1)$ for the rest of the proof.
By Part~2, the function $\cav u$ is uniformly continuous.
Hence,
there is $\delta > 0$ such that $|(\cav u)(p) - (\cav u)(p')| \leq \ep$ for every $p,p' \in \Delta(K)$ that satisfy $\|p - p'\|_\infty \leq \delta$. 
Assume w.l.o.g.~that $\delta < \ep^2$.
Then, $|(\cav u)(p) - (\cav u)(\xi_\delta(p))| \leq \ep$ for every $p \in \Delta(K)$.
%Fix $\delta < \ep^2$.

Set
\[ t_0 := 0, \ \ \ p_0 := p. \]
Consider now block $n$,
which starts
at stage $t_n$
with public history $h_{t_n}$
and initial belief $p_{t_n}(h_{t_n})$.

If $t_n = \theta$,
that is, if $p_{t_n}(h_{t_n})$ is close to the boundary of $\Delta(K)$,
we stop defining the blocks.

If $t_n < \theta$,
let $t_{n+1}$ be the first stage $t$ larger than $t_n$ such that at least one of the following conditions holds:
\begin{itemize}
\item
$t = \theta$.
\item
$\| p_{t_n}(h_{t_n}) - p_{t}(h_t)\|_1 \geq \delta \cdot \ep$.
\end{itemize}
If there is no integer $t$ larger than $t_n$ that satisfies any of these conditions, we set $t_{n+1} := \infty$.
In words, block $n$ starts at the history $h_{t_n}$ with a belief $p_{t_n}(h_{t_n})$.
The block continues as long as the posterior belief is far from the boundary of $\Delta(K)$ but close to $p_{t_n}(h_{t_n})$.

\bigskip
\noindent\textbf{Part~9: Definition of $\tau^*$.}

The strategy $\tau^*$ is defined in terms of the blocks that were described in Part~8.
For each $n$,
\begin{itemize}
\item 
If 
$t_n < \theta$,
then throughout block~$n$,
the strategy $\tau^*$ follows a subgame $\ep$-optimal
%minmax 
 strategy of Player~II in the non-revealing game $\Gamma^{\rm NR}(p_{t_n},f_{h_{t_n}})$.
\item
If 
$t_n=\theta$,
then 
$K_\delta(p_{t_n}(h_{t_n})) \subsetneq K$.
By the induction hypothesis,
$\underline v(\xi_\delta(\pi(p_{t_n},\sigma_{h_{t_n}},h_{t_n}))) = (\cav u)(\xi_\delta(\pi(p_{t_n},\sigma_{h_{t_n}},h_{t_n})))$.
Throughout block~$n$, the strategy $\tau^*$ follows 
a strategy that,
when facing $\sigma_{h_{t_n}}$, 
ensures that the expected payoff is at most 
%minimizes the payoff in 
$(\cav u)(\xi_\delta(\pi(p_{t_n},\sigma_{h_{t_n}},h_{t_n}))) + \ep$.
\end{itemize}

\bigskip

In the rest of the proof we show that the expected payoff under $(\sigma,\tau^*)$ is not much higher than $(\cav u)(p)$.

\bigskip

The next part implies that if the block is infinite with high probability, then the state
of nature does not significantly affect the public
play along the block.

\bigskip
\noindent\textbf{Part~10:
$| \prob_{k_1,\sigma(k_1), \tau^*}(A \mid h_T) - \prob_{k_1,\sigma(k_2), \tau^*}(A \mid h_T) | < 3\ep$,
for every $n \in \dN$,
every public history $h_T$ of length $T$ such%
\footnote{The inequality $t_n \leq T < t_{n+1}$ means that block $n$ does not end at the public history $h_T$.
In this inequality, the quantities $t_n$ and $t_{n+1}$ are in fact $t_n(h^\infty)$ and
$t_{n+1}(h^\infty)$ for any play $h^\infty \in H^\infty$ that extends $h_T$.}
that 
(i) $t_n \leq T < t_{n+1}$
and 
(ii) $\prob_{p_T(h_T),\sigma_{h_T},\tau^*_{h_T}}(t_{n+1} = \infty) \geq 1-\ep$,
every $t \geq T$,
every $k_1,k_2 \in K$,
and
every event $A \in \calF^{\rm II}_t$.}

Since $h_t$ extends $h_T$, by Bayes' rule we have
for every state of nature $k \in K$,

\begin{align*}
%\label{equ:1}
\begin{split}
p_t(h_t)(k)
&=  \prob_{p_{t_n}(h_{t_n}),\sigma_{h_{t_n}},\tau^*_{h_{t_n}}}(k^*=k \mid h_t)\\
&= \frac{ \prob_{k,\sigma_{h_{t_n}}(k),\tau^*_{h_{t_n}}}(h_t \mid h_T) \cdot  \prob_{p_{t_n}(h_{t_n}),\sigma_{h_{t_n}},\tau^*_{h_{t_n}}}(k^*=k \mid h_T)}{ \prob_{p_{t_n}(h_{t_n}),\sigma_{h_{t_n}},\tau^*_{h_{t_n}}}(h_t \mid h_T)}\\
&= \frac{ \prob_{k,\sigma_{h_{t_n}}(k),\tau^*_{h_{t_n}}}(h_t \mid h_T) \cdot p_T(h_T)(k)}{ \prob_{p_{t_n}(h_{t_n}),\sigma_{h_{t_n}},\tau^*_{h_{t_n}}}(h_t \mid h_T)}.
\end{split}
\end{align*}
Therefore, for every $k \in K$,
\begin{equation}
\label{equ:2}
 \prob_{k,\sigma_{h_{t_n}}(k),\tau^*_{h_{t_n}}}(h_t \mid h_T)
= \frac{p_t(h_t)(k) \cdot  \prob_{p_{t_n}(h_{t_n}),\sigma_{h_{t_n}},\tau^*_{h_{t_n}}}(h_t \mid h_T)}{p_T(h_T)(k)}.
\end{equation}

Since $A \in \calF^{\rm II}_t$,
there is a finite set $Z$ of public histories of length $t$ that extend $h_T$, such that
\[ A = \bigcup_{h_t \in Z} C(h_t).\]
Divide the set $Z$ into two subsets, $Z_0$ and $Z_1$,
according to whether along $h_t$, stage $t_{n+1}$ occurs 
after stage $t$.
Formally,
$Z_1$ contains all public histories $h_t \in Z$ such that, for some extension $h^\infty \in (I \times J)^\dN$ of $h_t$, we have $t_{n+1}(h^\infty) > t$,
and $Z_0 = Z \setminus Z_1$.

%$Z_0$ and $Z_1$,
%according to whether along $h_t$, stage $t_{n+1}$ occurs before stage $t$.
%Formally,
%$Z_1$ contains all public histories $h_t \in Z$ such that, for some extension $h^\infty \in (I \times J)^\dN$ of $h_t$, we have $t_{n+1}(h^\infty) > t$,
%and $Z_0 = Z \setminus Z_1$.

For any two distinct states of nature $k_1,k_2 \in K$,
%and any event $A \in \calF_t$ that is contained in $C(h_T)$,
\begin{align}
\nonumber
 \prob_{k_1,\sigma(k_1),\tau^*}&(A \mid h_T)
-  \prob_{k_1,\sigma(k_2),\tau^*}(A \mid h_T)\\
&=
\label{equ:5}
 \prob_{k_1,\sigma(k_1),\tau^*}(A \mid h_T)
-  \prob_{k_2,\sigma(k_2),\tau^*}(A \mid h_T)\\
\nonumber
&= \sum_{h_t \in Z}  \prob_{k_1,\sigma(k_1),\tau^*}(h_t \mid h_T)
- \sum_{h_t \in Z}  \prob_{k_2,\sigma(k_2),\tau^*}(h_t \mid h_T)\\
\nonumber
&= \sum_{h_t \in Z_1}  \prob_{k_1,\sigma(k_1),\tau^*}(h_t \mid h_T)
- \sum_{h_t \in Z_1}  \prob_{k_2,\sigma(k_2),\tau^*}(h_t \mid h_T)\\
\nonumber
&\hspace{1cm} + \sum_{h_t \in Z_0}  \prob_{k_1,\sigma(k_1),\tau^*}(h_t \mid h_T)
- \sum_{h_t \in Z_0}  \prob_{k_2,\sigma(k_2),\tau^*}(h_t \mid h_T)\\
\label{equ:11}
&\leq
\sum_{h_t \in Z_1}  \prob_{k_1,\sigma(k_1),\tau^*}(h_t \mid h_T)
- \sum_{h_t \in Z_1}  \prob_{k_2,\sigma(k_2),\tau^*}(h_t \mid h_T) + \ep\\
\label{equ:3}
&= \sum_{h_t \in Z_1}
 \prob_{p_{t_n}(h_{t_n}),\sigma_{h_{t_n}},\tau^*_{h_{t_n}}}(h_t \mid h_T)
\left(
\frac{p_t(h_t)(k_1)}{p_T(h_T)(k_1)} -
\frac{p_t(h_t)(k_2)}{p_T(h_T)(k_2)}
\right) + \ep\\
\label{equ:6}
&\leq
2|K|\ep + \ep \leq (2|K|+1)\ep,
\end{align}
where 
Eq.~\eqref{equ:5} holds since the posterior belief does not depend on the initial state
and since $A \in \calF^{\rm II}_t$,
Eq.~\eqref{equ:11} follows from the assumption that $\prob_{p_T,\sigma,\tau^*}(t_{n+1} = \infty \mid h_T) \geq 1-\ep$,
Eq.~\eqref{equ:3} follows from Eq.~\eqref{equ:2},
and
Eq.~\eqref{equ:6} holds since on
$Z_1$ we have
$|p_t(h_t)(k_1) - p_T(h_T)(k_1)| < \delta \cdot \ep$,
$|p_t(h_t)(k_2) - p_T(h_T)(k_2)| < \delta \cdot \ep$,
$p_T(h_T)(k_1) > \delta/|K|$, 
$p_T(h_T)(k_2) > \delta/|K|$,
and since $\delta < \ep^2$.

\bigskip

Since $\calF^{\rm II} = \sigma(\bigcup_{t \in \dN} \calF^{\rm II}_t)$,
the inequality proven in Part~10 holds for every event $A \in \calF^{\rm II}$  
that is contained in $\{t < \theta\}$,
and not only for events in a basis of $\calF^{\rm II}$.
That is,
for every $n \in \dN$,
every public history $h_T$ of length $T$ such that 
(i) $t_n \leq T < t_{n+1}$,
(ii) $ \prob_{p_T(h_T),\sigma_{h_T},\tau^*_{h_T}}(t_{n+1} = \infty) \geq 1-\ep$,
every event $A \in \calF^{\rm II}$, 
and every $k_1,k_2 \in K_\delta(p_T(h_T))$,
%every $n \in \dN$, every $T$ such that $t_n \leq T < t_{n+1}$,
%and every public history $h_t$ such that $p_t(h_t)$ is $\ep$-good
%and $ \prob_{p_T(h_T),\sigma_{h_T},\tau^*_{h_T}}(t_{n+1} = \infty) \geq 1-\ep$,
\begin{equation}
\label{equ:19}
\left|  \prob_{k_1,\sigma(k_1), \tau^*}(A \mid h_T) -  \prob_{k_1,\sigma(k_2), \tau^*}(A \mid h_T) \right| < (2|K|+1)\ep.
\end{equation}
Inequality~\eqref{equ:19} means that 
when the state of nature is $k_1$,
whether Player~I follows $\sigma(k_1)$ or $\sigma(k_2)$,
the two induced distributions on public plays that extend $h_T$ are close.
In other words, once the posterior process is required to be almost constant,
a condition that is captured by the assumption that $\prob_{p_T(h_T),\sigma_{h_T},\tau^*_{h_T}}(t_{n+1} = \infty) \geq 1-\ep$,
the knowledge of the state by Player~I does not affect much the distribution on future plays.

\bigskip

We use this observation to bound from above the expected payoff under $(\sigma_{h_T},\tau^*_{h_T})$.
This claim is not trivial, because in block $n$, the strategy $\tau^*$ is good when Player~I does not use the information she
possesses on the state,
but in fact Player~I does possess such information and may use it in $\sigma$.
Nevertheless,
it will follow from Part~10
that the fact that the posterior belief along block $n$ does not change by much
implies that $\sigma$ hardly uses any information on the state. In other words,
 $\sigma$ is not significantly different from a non-revealing strategy.

\bigskip
\noindent\textbf{Part~11:
$\E_{p_T(h_T),\sigma_{h_T},\tau^*_{h_T}}[f_{h_T}] \leq u(p_T(h_T)) + (2(2|K|+1)+1)\ep \leq u(p_T(h_T)) + 7|K|\ep$,
for every $n \in \dN$
and every public history $h_T$ of length $T$ such that 
(i) $t_n \leq T < t_{n+1}$
and 
(ii) $ \prob_{p_T(h_T),\sigma_{h_T},\tau^*_{h_T}}(t_{n+1} = \infty) \geq 1-\ep$.}

Fix $k_0 \in K$,
and define the strategy $\sigma^{\rm NR}$ in $\Gamma^{\rm NR}(p_T(h_T),h_T)$ by
\[ \sigma^{\rm NR} := \sigma_{h_T}(k_0). \]
That is, $\sigma^{\rm NR}$ follows $\sigma_{h_T}(k_0)$ whatever the state of nature is.
Since payoffs are between $0$ and $1$,
Eq.~\eqref{equ:19} implies that for every $k \in K_\delta(p_T(h_T))$,

\begin{align*}
&\left|\E_{k,\sigma_{h_T}(k_0),\tau^*_{h_T}}[f_{h_T}] - \E_{k,\sigma_{h_T}(k),\tau^*_{h_T}}[f_{h_T}] \right|\\ 
&=
\left| \int f_{h_T}(h^\infty) \rmd \prob_{\sigma_{h_T}(k_0),\tau^*_{h_T}}(h^\infty)
- 
\int f_{h_T}(h^\infty) \rmd \prob_{k,\sigma_{h_T}(k),\tau^*_{h_T}}(h^\infty)
\right|\\
&\leq
\int 
\left|\rmd \prob_{k,\sigma_{h_T}(k_0),\tau^*_{h_T}}(h^\infty)
-
\rmd \prob_{k,\sigma_{h_T}(k),\tau^*_{h_T}}(h^\infty)
\right|\\
&=
\left| 
\prob_{k,\sigma_{h_T}(k_0),\tau^*_{h_T}}(D)
-
\prob_{k,\sigma_{h_T}(k),\tau^*_{h_T}}(D)
\right|
+
\left| 
\prob_{k,\sigma_{h_T}(k_0),\tau^*_{h_T}}(D^{\rm c})
-
\prob_{k,\sigma_{h_T}(k),\tau^*_{h_T}}(D^{\rm c})
\right|
\\
&=
\left| 
\prob_{k,\sigma(k_0),\tau^*}(D \mid h_T)
-
\prob_{k,\sigma(k),\tau^*}(D \mid h_T)
\right|
+
\left| 
\prob_{k,\sigma(k_0),\tau^*}(D^{\rm c} \mid h_T)
-
\prob_{k,\sigma(k),\tau^*}(D^{\rm c} \mid h_T)
\right|
\\
&< 2(2|K|+1)\ep,
\end{align*}
where $D$ maximizes the difference
$\prob_{k,\sigma_{h_T}(k_0),\tau^*_{h_T}}(D)
-
\prob_{k,\sigma_{h_T}(k),\tau^*_{h_T}}(D)
$
among all sets in $\calF^{\rm II}$.

Since $\sigma^{\rm NR}$ follows $\sigma_{h_T}(k_0)$ whatever the state of nature is, this implies that

\begin{equation}
\label{equ:diff:small}
\left|\E_{p_T(h_T),\sigma^{\rm NR},\tau^*_{h_T}}[f_{h_T}] - \E_{p_T(h_T),\sigma_{h_T},\tau^*_{h_T}}[f_{h_T}] \right| < 2(2|K|+1)\ep.
\end{equation}
Since $\tau^*$ is a subgame $\ep$-minmax strategy in $\Gamma(p_{t_n}(h_{t_n}),f_{h_{t_n}})$,
the strategy
$\tau^*_{h_T}$ is $\ep$-minmax in $\Gamma(p_T(h_T),h_T)$, and hence
Eqs.~\eqref{equ:diff:small} and~\eqref{equ:value:independent} imply the desired result.

\bigskip
\noindent\textbf{Part~12}: The end of the proof.

For every $t \in \dN$, 
denote by $H^t_0$ the set of all plays $h^\infty$ such that 
either (a) $t \geq \theta(h^\infty)$,
or 
(b) $t \geq t_n$ and 
$\prob_{p_t,\sigma,\tau^*}(t_{n+1} = \infty \mid h_t) \geq 1-\ep$.
Since the posterior belief is a martingale,
it converges $\prob_{p,\sigma,\tau^*}$-a.s.
Hence, 
$\prob_{p,\sigma,\tau^*}(H^t_0) \geq 1-\ep$
for every $t$ sufficiently large.
Fix such $t$.
On the set $H^t_0 \cap \{t \geq \theta\}$
we have by the induction hypothesis,
\begin{equation}
\label{equ:121}
\E_{p,\sigma,\tau^*}[f \mid h_\theta] \leq 
(\cav u)(\xi_\delta(p_\theta(h_\theta))) + \ep
\leq (\cav u)(p_\theta(h_\theta)) + 2\ep, 
\end{equation} 
where the second inequality holds since 
$\cav u$ is uniformly continuous.
On the set $H^t_0 \cap \{t < \theta\}$ we have by Part~11,
\begin{equation}
\label{equ:122}
\E_{p,\sigma,\tau^*}[f \mid h_t] \leq u(p_t) + 7|K|\ep
\leq (\cav u)(p_t) + 7|K|\ep. 
\end{equation} 
Therefore,
\begin{align*}
\E_{p,\sigma,\tau^*}[f] 
&=
\E_{p,\sigma,\tau^*}\bigl[ \E_{p,\sigma,\tau^*}[f \mid \calF^{\rm II}_{t \wedge \theta}] \bigr]\\
&\leq 
\E_{p,\sigma,\tau^*}\bigl[ (\cav u)(p_{t \wedge \theta})] \bigr]+7|K|\ep\\
&\leq \E[(\cav u)(p)] + 8|K|\ep,
\end{align*}
where the first inequality holds by Eqs.~\eqref{equ:121} and~\eqref{equ:122}
and Jensen's inequality.
The second inequality holds since $\prob_{p,\sigma,\tau^*}(H^t_0) \geq 1-\ep$ and since payoffs are bounded by $1$.
Since $\ep$ is arbitrary, the result follows.

%%%%%%%%%%%%%%%%%%%%%%%%%%%%%%%%%%%%%%%%%
\section{Examples}
\label{section:examples}

In this section we provide two examples,
which are variants of a game studied by Laraki (2000). The first example shows that tail measurability does not guarantee the existence of the value.
We will consider a game with alternating moves:
Player~I selects actions \emph{only} in even stages,
and Player~II selects actions \emph{only} in odd stages.
In particular, a play is a vector $(k,i_0,j_1,i_2,j_3,\dots)$,
where $k \in K$, $i_{2t} \in I$ and $j_{2t+1} \in J$ for every $t \in \dN$.
Taking the identity of the players into account,
a function $f$ defined on the set of plays is \emph{shift invariant}
%, and hence tail measurable,
if $f(k,i_0,j_1,i_2,j_3,\dots) = f(k,i_2,j_3,i_4,j_5,\dots)$
for every play $(k,i_0,j_1,i_2,j_3,\dots)$.
The examples provided before for tail-measurable functions are in fact shift invariant.
Any shift-invariant function is tail measurable, but the converse is not true.
For example, the function $f$ that is equal to 1 if Player~I plays a specific action $i_*$ in infinitely many periods $t$ such that $t = 0$ modulo 4, and 0 otherwise, is tail measurable but not shift invariant.

\subsection{Example 1}
\label{section:example1}

Suppose that $K = \{k_1,k_2\}$, $I = J = \{\ell,r\}$, and
the payoff function is given by %(see also Figure~\arabic{figurecounter}):
\[f(k,i_0,j_1,i_2,j_3\dots) =\begin{cases}
            -1,  &\textnormal{if }k=k_{1},\:j_{2t+1}=\ell\:\textnormal{infinitely often,} \\
            2,  &\textnormal{if } k=k_{2},\:j_{2t+1}=\ell\:\textnormal{infinitely often,} \\
            -2,  &\textnormal{if } k=k_{1},\:j_{2t+1}=\ell\:\textnormal{for finitely many stages}, \:i_{2t}=\ell\:\textnormal{infinitely often,} \\
            1,  &\textnormal{if } k=k_{2},\:j_{2t+1}=\ell\:\textnormal{for finitely many stages}, \:i_{2t}=\ell\:\textnormal{infinitely often,}\\
            0, & \textnormal{otherwise,}
        \end{cases}
    \]
which is shift invariant.

%\[
%\begin{array}{lllr}
%\hbox{State} & \hbox{Player~I} & \hbox{Player~II} & \hbox{Payoff}\\
%\hline
%k_1 & \hbox{does not matter} & \hbox{quits i.o.} & -1 \\
%k_1 & \hbox{quits i.o.} & \hbox{quits finitely many times} & -2 \\
%k_1 & \hbox{quits finitely many times}& \hbox{quits finitely many times} & 0\\
%k_2 & \hbox{does not matter} & \hbox{quits i.o.} & 2 \\
%k_2 & \hbox{quits i.o.} & \hbox{quits finitely many times} & 1 \\
%k_2 & \hbox{quits finitely many times}& \hbox{quits finitely many times} & 0\\
%\end{array}
%\]
%\centerline{Figure~\arabic{figurecounter}: The payoff function in the game in Section~\ref{section:example:no:value}.}
%\addtocounter{figurecounter}{1}

The value function $u$ of the non-revealing game
is (see Figure~\arabic{figurecounter}):
    \begin{equation}
        u(p) = \begin{cases}
                1-3p,  & \textnormal{if }0\le p\le\frac{1}{3}, \\
                0, & \textnormal{if } \frac{1}{3}\le p\le\frac{2}{3}, \\
                2-3p,  & \textnormal{if }\frac{2}{3}\le p\le1. \\
            \end{cases}
    \end{equation}

\centerline{\includegraphics[width=8cm]{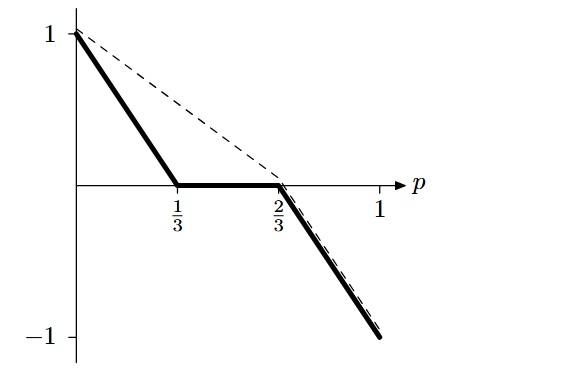}}
%\centerline{\includegraphics{figure.1}}
\centerline{Figure~\arabic{figurecounter}: The functions $u$ (dark) and $\cav u$ (dashed).}
\addtocounter{figurecounter}{1}

\bigskip

Indeed, when $p\in[0,\frac{1}{3}]$, we have $2-3p\geq1-3p\ge0$.
When Player~II always selects $r$, the payoff is at most $\max\{0,1-3p\}\leq1-3p$, hence $u(p)\leq1-3p$.
When Player~I always selects $\ell$, the payoff is at least $\min\{1-3p,2-3p\}\geq1-3p$, hence $u(p)\geq1-3p$.

When $\frac{1}{3}<p<\frac{2}{3}$, we have $1-3p<0<2-3p$.
When Player~II always selects $r$, the payoff is at most $\max\{0,1-3p\}\leq0$, hence $u(p)\leq0$.
When Player~I always selects $r$, the payoff is at least $\min\{0,2-3p\}\geq0$, hence $u(p)\geq0$.

The case $\frac{2}{3}\le p\le1$ is analogous to the case $p\in[0,\frac{1}{3}]$.

\bigskip

We will show that the game $\Gamma(\frac{1}{2})$ admits no value.
By Theorem~\ref{theorem:main},
$\underline v(\frac{1}{2}) = (\cav u)(\frac{1}{2}) = \frac{1}{4}$.
Hence, to show that $\Gamma(\frac{1}{2})$ admits no value,
it is sufficient to show that
$\overline v(\frac{1}{2}) \geq \frac{1}{2} > \frac{1}{4}$.

Fix $\ep>0$ and a strategy $\tau$ for Player~II.
Denote by $\sigma_r$ the strategy for Player~I in which she always selects $r$, independently of the state, the stage, and the players' past actions.
Let $E_1$ be the event that Player~II selects $\ell$ infinitely often.
There exist $t \in \dN$ and an event $\widehat E_1 \in \calF^{\rm II}_t$ such%
\footnote{Recall that $E_1 \triangle E_2 = (E_1 \setminus E_2) \cup (E_2 \setminus E_1)$, for any two sets (events) $E_1,E_2$.}
that $\prob_{p,\sigma_r,\tau}(E_1 \triangle \widehat E_1) < \ep$.
Thus, if $h_t \in \widehat E_1$,
then with high probability Player~II will select $\ell$ infinitely often,
and if $h_t \not\in \widehat E_1$,
then with high probability Player~II will select $\ell$ only finitely many times.

A strategy $\sigma$ of Player~I that guarantees $\frac{1}{2} - \ep$ against $\tau$ suggests itself:
\begin{itemize}
    \item If $k=k_1$, always play $r$.
    \item If $k=k_2$, play $r$ until stage $t$, and afterwards, on the event $\widehat E_1$ continue selecting $r$, and on the event $(\widehat E_1)^{\rm c}$ select $\ell$ in all subsequent stages.
\end{itemize}

\subsection{Example 2}
\label{section:example2}

The payoff function in 
the game presented in Section~\ref{section:example1}
is such that two plays that differ in a set of action pairs that has density 0 yield different payoffs.
This does not happen for the long-run average payoff function.
One might expect that the value would exist when the payoff function is shift invariant
and has the property that every two plays that differ in a set of action pairs with density 0 yield the same payoff.
However, this is not the case.
To see this, we consider a variant of the game in Section~\ref{section:example1}.
In the context of that example, define
\begin{align*}
d_1(i_0,i_2,\dots) &:= \limsup_{t \to \infty} \frac{1}{t} \sum_{m=0}^{t-1} \mathbf{1}_{\{i_{2m} = \ell\}}, \\
d_2(j_1,j_3,\dots) &:= \limsup_{t \to \infty} \frac{1}{t} \sum_{m=0}^{t-1} \mathbf{1}_{\{j_{2m+1} = \ell\}}.
\end{align*}
These are the densities of the set of stages in which each player selects $\ell$.
Let $f$ be the following payoff function:
\[f(k,i_0,j_1,i_2,j_3,\dots) =\begin{cases}
            -1,  &\textnormal{if }k=k_{1},\: d_2(j_1,j_3,\dots) > 0.1, \\
            2,  &\textnormal{if } k=k_{2},\: d_2(j_1,j_3,\dots) > 0.1, \\
            -2,  &\textnormal{if } k=k_{1},\:d_2(j_1,j_3,\dots) \leq 0.1, \: d_1(i_0,i_2,\dots) > 0.1,\\
            1,  &\textnormal{if } k=k_{2},\:d_2(j_1,j_3,\dots) \leq 0.1, \: d_1(i_0,i_2,\dots) > 0.1,\\
            0, & \textnormal{otherwise.}
        \end{cases}
    \]

The arguments in Section~\ref{section:example1} show that the minmax and maxmin values of this game coincide with those of the game in Section~\ref{section:example1}.

%%%%%%%%%%%%%%%%%%%%%%%%%%%%%
\section{Open Problems}
\label{section:discussion}

In this section we list some open problems that motivated this study.

The concavification of the value function of the non-revealing game plays an important role in this paper.
For the long-run average payoff,
AM provided a simple characterization of this function
that involves only the value of proper one-shot games.
It would be interesting to find other general classes of payoff functions for which concavification of the value function of the non-revealing game can be explicitly or approximately characterized.

Repeated games with incomplete information on one side admit a value 
under the long-run average payoff, 
and do not admit a value when the payoff function is shift invariant. 
A natural question concerns the properties of the payoff function that guarantee the existence of the value.

Another open problem relates to the characterization of the maxmin value.
For some payoff functions, the maxmin value coincides with the concavification of the value function of the non-revealing games.
In addition to games with long-run average payoffs,
these include, e.g., 
certain variations of the Big Match game, though not all of them (Sorin, 1984, 1985), as well as deterministic stopping games (Laraki, 2000).
Under what conditions on the payoff function is the maxmin value equal to the concavification of the value function of the non-revealing game?

%\blue{instead of :}
%payoff functions that arise from 
%certain types of Big Match games (Sorin, 1985)
%and deterministic stopping games (Laraki, 2000).
%For other payoff functions, this is not the case, 
%e.g., payoff functions that arise from
%other types of Big Match games (Sorin, 1984).
%One then seeks the properties of the payoff function that guarantee that the maxmin value coincides with the concavification of the value function of the non-revealing game.

Finally, 
one would like to characterize the maxmin value when 
it does not coincide with the concavification of the value function of the non-revealing game,
as well as the minmax value of the game.

%\blue{Necessary?  :}We hope that our study will pave the way to answering these open problems.

\end{document}